\documentclass[11pt,a4paper]{article}
\usepackage{amssymb,amsmath}
\usepackage[abbrev]{amsrefs}
\usepackage{latexsym}
\usepackage{graphicx,color}
\usepackage{times,theorem,latexsym,color}

\newcommand{\BOX}{\ensuremath\Box}

\newtheorem{theorem}{Theorem }[section]

{\theorembodyfont{\rmfamily}}
{\theorembodyfont{\rmfamily}}
{\theorembodyfont{\rmfamily}}
\newtheorem{lemma}[theorem]{Lemma}

{\theorembodyfont{\rmfamily}\newtheorem{remark}[theorem]{Remark}}
{\theorembodyfont{\rmfamily}}

\newcommand{\R}{\mathbb{R}}

\newcommand{\dd}{\,{\rm d}}

{\hspace*{.1pt}\hspace*{\fill}\BOX\vskip\baselineskip}

{\vskip\baselineskip\noindent\textbf{Proof of {#1}:}}%
{\hspace*{.1pt}\hspace*{\fill}\BOX\vskip\baselineskip}
{\vskip\baselineskip\noindent\textbf{Proof of Theorem \protect\ref{#1}:}}%
{\hspace*{.1pt}\hspace*{\fill}\BOX\vskip\baselineskip}
{\vskip\baselineskip\noindent\textbf{Proof of Theorems \protect\ref{#1} --
\protect\ref{#2}:}}%
{\hspace*{.1pt}\hspace*{\fill}\BOX\vskip\baselineskip}

\begin{document}

\title{On the two-dimensional steady Navier-Stokes equations related to flows around a rotating obstacle}
\author{Mitsuo Higaki\thanks{Department of Mathematics, Kyoto University; \texttt{mhigaki@math.kyoto-u.ac.jp}} \and Yasunori Maekawa\thanks{Department of Mathematics, Kyoto University; \texttt{maekawa@math.kyoto-u.ac.jp}} \and Yuu Nakahara\thanks{Mathematical Institute, Tohoku University; \texttt{yuu.nakahara.t3@dc.tohoku.ac.jp}}}

\date{}

\maketitle

\section{Introduction}
Let $\mathcal B$ be a rigid body immersed in a viscous incompressible fluid that fills the whole space. Assume that the body rotates with a constant angular velocity $a\in \R\setminus\{0\}$ and the exterior of $\mathcal B(t)$ is described as $\Omega(t) \subset \R^2$. The time dependent domain $\Omega(t)$ is defined as 
	\begin{align}\label{rota}
		\begin{split}
			\Omega(t) & \, = \, \big \{ y\in \R^2~|~ y = O (a t) x\,, x\in \Omega \big \}\,,\\
			O (a t ) & \, = \, 
				\begin{pmatrix}
					\cos a t & -\sin a t\\
					\sin a t & \cos a t
				\end{pmatrix}\,,
		\end{split}
	\end{align}

\noindent where a given exterior domain $\Omega(0) = \Omega \subset\R^2$ has a smooth boundary $\partial \Omega$.
The flow around the rotating body is described by the following Navier-Stokes equations:
	\begin{equation}\label{NS}
	  \left\{
		\begin{aligned}
			\partial_t v -\Delta v + v\cdot \nabla v + \nabla q & \,=\, g \,, \qquad t>0\,,~y \in \Omega (t)\,, \\
			{\rm div}\, v & \,=\, 0 \,, \qquad  t>0\,,~ y \in \Omega (t)\,. \\
		\end{aligned}\right.
	\end{equation}
	
\noindent Here $v=v(y,t) = (v_1(y,t), v_2 (y,t))^\top$ and $q=q(y,t)$ are respectively unknown velocity field and pressure field, and $g = g(y,t) = (g_1(y, t), g_2(y,t))^\top$ is a given external force. 
We use the standard notation for derivatives: $\partial_t = \frac{\partial}{\partial t}$, $\partial_j = \frac{\partial}{\partial x_j}$, $\Delta = \sum_{j=1}^2 \partial^2_j$, ${\rm div}\, v = \sum_{j=1}^2 \partial_j v_j$, $v\cdot \nabla v = \sum_{j=1}^2 v_j \partial_j v$, while $x^{\bot}=(-x_2,x_1)^\top$ denotes the vector which is perpendicular to $x =(x_1,x_2)^\top$. To get rid of the difficulty due to the time-dependence of the domain we take the reference frame by making change of variables for $t\geq 0$ and $x\in \Omega$
	\begin{align*}
		y \, = \, O (a t) x\,, \quad u (x,t) &\, = \, O (a t)^\top v (y,t)\,, \quad p (x,t) \, = \, q (y,t)\,, \\  f (x,t) &\, = \, O (a t)^\top g (y,t)\,.
	\end{align*} 

\noindent 
Then \eqref{NS} is equivalent to the equations:
	\begin{equation}\label{NSNS}
		\left\{
			\begin{aligned}
				\partial_t u -\Delta u - a ( x^\bot \cdot \nabla u - u^\bot ) + \nabla p & \,=\, -u\cdot \nabla u + f \,, \quad & t>0\,,~x \in \Omega \,, \\
				{\rm div}\, u & \,=\, 0\,, & t>0\,,~ x \in \Omega \,. \\
			\end{aligned}\right.
	\end{equation}

\noindent  In order to understand the structure of solutions at spatial infinity it is important to study this system in $\R^2$. The effect of the boundary is expressed as a force in this case. Motivated by this observation, as a model problem, in this paper we study the above nonlinear system in $\R^2$ and in the steady case. Thus, assuming that $f$ is independent of $t$, we are interested in the following system:
\begin{equation}\tag{NS$_a$}\label{NS_a}
		\left\{
			\begin{aligned}
			 	-\Delta u - a ( x^\bot \cdot \nabla u - u^\bot ) + \nabla p & \,=\, 
			 	- u\cdot\nabla u  + f \,, \quad\  &x \in \R^2\,,\\
				{\rm div}\, u \,& \,=\, 0\,,  &x \in \R^2\,. \\
			\end{aligned}
		\right.
	\end{equation}
Our aim is to show the existence and the asymptotic behavior of the solution to \eqref{NS_a}. For this purpose we first consider the linearized problem
	\begin{equation}\tag{S$_a$}\label{S_a}
		\left\{
			\begin{aligned}
				 -\Delta u - a ( x^\bot \cdot \nabla u - u^\bot ) + \nabla p & \,=\, 
				 f \,, \quad\  &x \in \R^2\,,\\
				 {\rm div}\, u \,& \,=\, 0\,,  &x \in \R^2\,. \\
			\end{aligned}
		\right.
	\end{equation}
We will show that there exists a unique solution to \eqref{S_a} such that the leading term of the asymptotic behavior of the flow at infinity is the rotational profile $c \frac{x^\bot}{|x|^2}$ whose coefficient $c$ is determined by the external force $f$.

Before stating the main theorem, 
let us recall some known results on the mathematical  analysis of flows around a rotating obstacle.

So far the mathematical results on this topic have been obtained mainly for the three-dimensional problem,
as listed below.
For the nonstationary problem the existence of global weak solutions is proved by Borchers \cite{Bo}, and the unique existence of time-local regular solutions is shown by Hishida \cite{H1} and Geissert, Heck, and Hieber \cite{GHH}, while the global strong solutions for small data are obtained by Galdi and Silvestre \cite{GSi}. The spectrum of the linear operator  related to this problem is studied by Farwig and Neustupa \cite{FN}; see also the linear analysis by Hishida \cite{H2}. The existence of stationary solutions to the associated system is proved in \cite{Bo}, Silvestre \cite{Si},  Galdi \cite{G1}, and Farwig and Hishida \cite{FH0}. In particular, in \cite{G1} the stationary flows with the decay order  $O(|x|^{-1})$ are obtained, while the work of \cite{FH0} is based on the weak $L^{3}$ framework, which is another natural scale-critical space for the three-dimensional Navier-Stokes equations. 
In $3$D case the asymptotic profiles of these stationary flows at spatial infinity are studied by Farwig and Hishida \cite{FH1,FH2} and Farwig, Galdi, and Kyed \cite{FGK}, where it is proved that the asymptotic profiles are described by the Landau solutions, stationary self-similar solutions to the Navier-Stokes equations in 
$\R^3\setminus\{0\}$. 
It is worthwhile to mention that, also in the two-dimensional case,
the asymptotic profile is given by the stationary self-similar solution $c \frac{x^\bot}{|x|^2}$.
The stability of the above stationary solutions has been well studied in the three-dimensional case;
The global $L^2$ stability is proved in \cite{GSi}, 
and the local $L^3$ stability is obtained by Hishida and Shibata \cite{HShi}.

All results mentioned above are considered in the three-dimensional case, 
while only a few results are known so far for the flow around a rotating obstacle in the two-dimensional case. 
An important progress has been made by Hishida \cite{H3}, 
where the asymptotic behavior of the two-dimensional stationary Stokes flow around a rotating obstacle is investigated in details.
Recently, the nonlinear problem \eqref{NS} is analyzed in \cite{HMN}, and the existence of the unique solution decaying as $O(|x|^{-1})$ is proved for sufficiently small $a$ and $f$ when the external force $f$ is of divergence form $f={\rm div}\, F$ and $F$ has a scale critical decay. Moreover, the leading profile at spatial infinity is shown as $C \frac{x^\bot}{|x|^2}$ under the additional decay condition on $F$ such as $F=O(|x|^{-2-r})$, $r>0$. 

Since we consider the problem in $\R^2$ in this paper, by virtue of the absence of the physical boundary,
we can show the existence of solutions to \eqref{NS_a} without assuming the smallness of the angular velocity $a$. To state our result let us introduce the function space. 
For a fixed number $s\geq 0$ the weighted $L^\infty$ space $L^\infty_s (\R^2)$ is defined as
	\begin{align*}\label{def.L^infty_s}
		L^\infty_s (\R^2) & \,=\, \big \{ f\in L^\infty (\R^2)~|~ (1+|x|)^s f \in L^\infty (\R^2) \big \}\,.
	\end{align*}
The space is a Banach space equipped with the natural norm 
\begin{equation*}
\|f\|_{L^\infty_s} \,=\, {\rm ess.sup}_{x\in\R^2} (1+|x|)^s |f(x)|\,.
\end{equation*}
The first result of this paper is stated as follows. 
\begin{theorem}\label{thm.main} 
Let $a \in \R \setminus \{ 0 \}$ and $r\in[0,1). $ Assume that $f \in L^\infty_{3+r}(\R^2)^2$. 
Then there exists a unique $(u,p) \in L^\infty_1 (\R^2)^2 \times L^\infty (\R^2)$ such that:
	\begin{enumerate}
		\item The couple $(u,p)$ satisfies \eqref{S_a} in the sense of distributions.
		\item The velocity $u$ belongs to $L^\infty_1(\R^2)^2$ and satisfies
			\begin{align}\label{1.1}
				u(x) \,=\,  \int_{|y|<\frac{|x|}2} y^\bot \cdot f(y) \dd y \ \frac{ x^{\bot} }{4\pi |x|^2} + \mathcal{R}[f](x)\,, 
			\end{align}
		with
\begin{align}\label{1.2}
				|\mathcal{R}[f](x)| \leq \frac{C}{(1+|x|)^{1+r}}\big (\frac{1}{1-r} + \frac{1}{|a|^\frac{1+r}{2}} \big )   ||f||_{L^\infty_{3+r}}\,,
\end{align}
        and in particular, it follows that 
\begin{align}\label{1.3}
\| \mathcal{R}[f]\|_{L^\infty_{1+r}} & \leq C \big (\frac{1}{1-r}+\frac{1}{|a|^\frac{1+r}{2}} \big ) \| f\|_{L^\infty_{3+r}}
\end{align}
with a numerical constant $C$.
   		\item The pressure $p$ is given by
   			\begin{align}\label{bangou}
   				p(x) \,=\, \frac{1}{2\pi} \int_{\R^2} \frac{x-y}{|x-y|^2} \cdot f(y) \dd y\,.
   			\end{align}
	\end{enumerate}
\end{theorem}

\begin{remark} (1) The representation \eqref{bangou} leads to the regularity of the pressure such as  $p\in L^\infty_1 (\R^2)$. The solution $(u,p)\in L^\infty_1(\R^2)^2 \times L^\infty_1(\R^2)$ satisfying \eqref{S_a} in the sense of distributions is unique by virtue of the uniqueness result  in Hishida \cite[Lemma 3.5]{H3}. 

\noindent (2) In \cite{H3} the result of Theorem \ref{thm.main} is firstly established under the conditions on $f$ such as $f \in L^1(\R^2)^2 \cap L^\infty(\R^2)^2,\ x^\bot \cdot f \in L^1(\R^2)^2$, and $f(x)=O(|x|^{-3}(\log{|x|})^{-1})$ as $|x|\rightarrow \infty$. Our result improves his result, and in particular, the critical case $f = O(|x|^{-3})$ is treated. We note that, in the case $r=0$, the integral in the right-hand side of \eqref{1.1} does not converge in general when $|x|\rightarrow \infty$.
\end{remark}

To study the nonlinear problem it is reasonable to consider the linear problem when the external force $f$ is given by $f={\rm div}\,F$ with $F(x) = O(|x|^{-2})$ in view of the structure of the nonlinear term $u\cdot \nabla u = {\rm div}\, (u \otimes u)$. Here the matrix $(u_i v_j)_{1\leq i,j \leq 2}$ is written as $u \otimes v$. The following result is essentially obtained in \cite{HMN}.

\begin{theorem}[ {\rm \cite[Theorem 3.1(ii)]{HMN}}]\label{lem.F} Let $a \in \R \setminus \{ 0 \}$, $r\in[0,1)$, and $q\in (1,\infty)$. Assume that  $f\in L^2(\R^2)^2$ is of divergence form $f={\rm div}\, F=(\partial_1 F_{11} + \partial_2 F_{12}, \partial_1 F_{21} + \partial_2 F_{22})^\top$ with $F=(F_{ij})_{1\leq i,j\leq 2} \in L^\infty_{2+r}(\R^2)^{2\times2} $. Then there exists a unique $(u,p) \in L_1^\infty(\R^2)^2 \times L^q (\R^2)$ such that $(u,p)$ satisfies \eqref{S_a} in the sense of distributions, and $u$ satisfies 
\begin{align}\label{thelemma2}
u(x) \,=\, \int_{|y|<\frac{|x|}2} (F_{12}(y) - F_{21}(y)) \dd y \frac{x^\bot}{4\pi |x|^2} + \mathcal{R}[f](x)\,,
\end{align}
and $\mathcal{R}[f]$ satisfies
\begin{align}\label{thelemma3}
\begin{split}
|\mathcal{R}[f](x)| & \leq C  \min \big\{ \frac{1}{|a| |x|^3}, \frac{1}{|x|} \big\}  \int_{|y|\leq \frac{|x|}{2}} |F(y) | \dd y \\
& \qquad + \frac{C}{(1+|x|)^{1+r}} \frac{1}{1-r} \|F\|_{L^\infty_{2+r}}\,.
\end{split}
\end{align}
Here $C$ is a numerical constant independent also of $a$ and $r$.
In particular, it follows that 
\begin{align}\label{1.4}
\| \mathcal{R}[f]\|_{L^\infty_{1+r}} \leq C \big (\frac{1}{1-r} + \frac{1 + \log |a|}{|a|^\frac{r}{2}} \big ) \| F \|_{L^\infty_{2+r}}
\end{align}
with a numerical constant $C$.
\end{theorem}

\begin{remark}{\rm (1) In fact, the statement of \cite[Theorem 3.1]{HMN} is a slightly different from Theorem \ref{lem.F} above. So we give a sketch of the proof of Theorem \ref{lem.F} in Section \ref{sec.pre}, based on the key pointwise asymptotic estimate of the fundamental solution, see Lemma \ref{lem.thm.linear.whole.1} below, which is due to \cite[Lemma 3.3]{HMN}.

\noindent (2) Estimate \eqref{1.4} is derived from \eqref{thelemma3}. Indeed, when $|x|\leq 1$ the first term in the right-hand side of \eqref{thelemma3} is estimated as $C|x|\|F\|_{L^\infty}$, while when $|x|\geq 1$ this term is estimated by dividing into two cases (i) $|a||x|^2\leq 1$ and (ii) $|a||x|^2\geq 1$. The factor $\log |a|$ in \eqref{1.4} is required only when $r$ is near $0$ in order to ensure that the constant $C$ is independent of $r$.
}
\end{remark}	

The linear results of Theorem \ref{thm.main} and \ref{lem.F} are applied to the nonlinear problem \eqref{NS_a}. The result for the nonlinear problem is stated as follows.
\begin{theorem}\label{thm.main2} Let $a \in \R \setminus \{ 0 \}$ and $r\in [0,1)$. 
Then there exists $\delta=\delta(a,r)>0$ such that, for any $f \in L^\infty_{3+r}(\R^2)^2$ satisfying $x^\bot \cdot f\in L^1 (\R^2)$ and 
	\begin{align}\label{delta}
	 		\| x^\bot \cdot f\|_{L^1} + \| f \|_{L^\infty_{3+r}} < \delta\,,
	 \end{align}
there exists a unique solution $(u,p)$ to \eqref{NS_a} such that
\begin{align}
\begin{split}
u(x) \,=\, \alpha U(x) + v(x)\,,
\end{split}
\end{align}
where 
	\begin{align}\label{alpha} 
			\alpha \,=\, \frac12 \int_{\R^2} y^{\bot} \cdot f(y)\ \dd y\,, \qquad U(x) \,=\, \frac1{2\pi} \frac{x^\bot}{|x|^2}&(1-e^{-\frac{|x|^2}4})\,,
	\end{align}
and
\begin{align}\label{1.5}
\| v \|_{L^\infty_{1+r}} & \leq C_r \big  (\frac{1}{1-r} + \frac{1}{|a|^\frac{1+r}{2}} \big) \big (\| x^\bot \cdot f \|_{L^1} + \| f \|_{L^\infty_{3+r}} \big )\,,
\end{align}
and the pressure $p$ is given by 
	\begin{align}\label{bangou2}
			p(x) \,=\, \nabla \cdot (-\Delta)^{-1} \nabla \cdot (u \otimes u) - \nabla \cdot (-\Delta)^{-1} f\,.
	\end{align}
Here the constant $C_r$ depends only on $r$.
\end{theorem}
\begin{remark}{\rm In Theorem \ref{thm.main2} the solution is constructed as the solution to the integral equation associated with \eqref{NS_a}, which is formulated based on the fundamental solution to the linearized problem \eqref{S_a}. The uniqueness is proved for this class of solutions.
}
\end{remark}

This paper is organized as follows. In Section \ref{sec.pre} we collect the estimates which reflect the effect of the rotation. Most of them are the abstractions from \cite{H3,HMN}. 
Theorem \ref{thm.main} is proved in Section \ref{sec.thm.main}. Finally, Theorem \ref{thm.main2} is proved in Section \ref{sec.nonlinear}.

\section{Preliminaries}\label{sec.pre}
Let us consider the linear problem in the whole plane for $a \in \R \setminus \{0\}$:
	\begin{equation}\tag{S$_a$}\label{a}
			- \Delta u -a (x^\bot \cdot \nabla u-u^\bot)+ \nabla p \,=\, f\,, \qquad {\rm div}\ u \,=\, 0\,,  \qquad \quad x \in  \R^2\,.
	\end{equation}
\noindent 
Let $q\in [1,\infty]$.
The couple $(u,p)\in L^\infty (\R^2)^2\times L^q(\R^2)$ is said to be a weak solution to \eqref{a} 
if (i) ${\rm div}\, u=0$ in the sense of distributions, and (ii) $(u,p)$ satisfies
	\begin{equation*}\label{def.weak.whole}
			\int_{\R^2}  u\cdot  T_{-a} \phi  \dd x  - \int_{\R^2} p \, {\rm div}\, \phi \dd x  \,=\,  \int_{\R^2} f \cdot \phi \dd x\,, \qquad {\rm for~all} \ \ \phi\in \mathcal{S}(\R^2)^2\,,
	\end{equation*}
where
\begin{align*}
	T_a u \,=\, -\Delta u -a(x^\bot \cdot \nabla u -u^\bot)\,.
\end{align*}
Let $\mathbb{I}=(\delta_{ij})_{1\leq i,j \leq 2}$ be the identity matrix.
The velocity part of the fundamental solution to \eqref{a} plays a central role throughout this paper, which is defined as 
\begin{equation}\label{gamma}
			\Gamma_{a}(x,y) \,=\,  \int_{0}^{\infty} O(a t)^{\top}K(O(a t)x-y,t) \dd t\,,
\end{equation}
where
\begin{align}\label{kernel}
			K(x,t) \, = \, G(x,t) \mathbb{I} + H(x,t)\,, \qquad  H(x,t) \, = \, \int_{t}^{\infty} \nabla^2 G(x,s) \dd s\,, 
\end{align}
and $G(x,t)$ is the two-dimensional Gauss kernel
\begin{equation*}
			G(x,t) \,=\, \frac{1}{4\pi t} e^{-\frac{|x|^2}{4t}}\,.
\end{equation*}
Similarly, the pressure part of the fundamental solution is defined as 
\begin{equation*}
			Q(x-y) \,=\, \frac1{2\pi}\log{|x-y|}\,,
\end{equation*}
for the following identity holds.
\begin{equation*}
			{\rm div}\, (x^\bot \cdot \nabla u-u^\bot )=x^\bot \cdot \nabla{\rm div}\, u \,=\, 0\,.
\end{equation*}
\begin{remark}
We can also write $H(x,t)$ in \eqref{kernel} as follows
\begin{align*}
H(x,t) \,=\, -\frac{(x \otimes x)}{|x|^2} G(x,t) + \bigg(
				\frac{x \otimes x}{|x|^2} - \frac{\mathbb I}{2}
			\bigg) \frac{1- e^{-\frac{|x|^2}{4t}}}{\pi |x|^2}\,.
\end{align*}	
\end{remark}
%
The next lemma is proved in \cite{H3,HMN}.
\begin{lemma}[ {\cite[Proposition 3.1]{H3}, \cite[Lemma 3.3]{HMN}}]\label{lem.thm.linear.whole.1} Set
\begin{equation}\label{def.Lxy}
L(x,y) \,=\, \frac{x^\bot \otimes y^\bot}{4\pi |x|^2}\,. 
\end{equation}
Then for $m=0,1$ the kernel $\Gamma_{a}(x,y)$ satisfies
\begin{equation}\label{est.lem.thm.linear.whole.1.1}
\begin{split}
& | \nabla_y ^m \big ( \Gamma_{a}(x,y)- L(x,y) \big ) | \\
& \le
C \bigg ( \delta_{0m} \min \big\{ \frac{1}{|a| |x|^2}, \frac{1}{|a|^{\frac{1}{2}} |x|} \big\}
+  |x|^{1-m} \min \big\{ \frac{1}{|a| |x|^3}, \frac{1}{|x|} \big\} 
+  \frac{|y|^{2-m}} {|x|^2}\bigg )\,,\\
&  \quad\quad\quad\quad\quad \quad\quad\quad\quad\quad\quad \quad\quad \quad\quad\quad \quad\quad\quad\quad\quad  {\rm for}  \quad |x|> 2 |y|\,.
\end{split}
\end{equation}
Here $\delta_{0m}$ is the Kronecker delta and $C$ is independent of $x$, $y$, and $a$. 
\end{lemma}

\begin{remark}\label{rem.lem.thm.linear.whole.1}{\rm (1) The asymptotic estimate like \eqref{est.lem.thm.linear.whole.1.1} is proved in \cite{H3} when $m=0$, and then the dependence on $|a|$ is improved by \cite{HMN} which is needed to solve the nonlinear problem. The detailed proof for the case $m=1$ of \eqref{est.lem.thm.linear.whole.1.1} is given by \cite{HMN}.
 
\noindent (2) Note that, when $|y|>2|x|$, since $\Gamma_a(x,y)=\Gamma_{-a}(y,x)^\top $ and $(y^\bot \otimes x^\bot)^\top=x^\bot \otimes y^\bot$ we have a similar estimate:
\begin{align}\label{reverce lem}
\begin{split}
& | \Gamma_{a}(x,y)- \frac{x^\bot \otimes y^\bot}{4\pi |y|^2} | \\
& \le
C \bigg (\min \big\{ \frac{1}{|a| |y|^2}, \frac{1}{|a|^{\frac{1}{2}} |y|} \big\}
+  |y| \min \big\{ \frac{1}{|a| |y|^3}, \frac{1}{|y|} \big\} 
+  \frac{|x|^{2}} {|y|^2}\bigg )\,,\\
&  \quad\quad\quad\quad\quad \quad\quad\quad\quad\quad\quad \quad\quad \quad\quad\quad \quad\quad\quad\quad\quad  {\rm for}  \quad |y|> 2|x|\,.
\end{split}
\end{align}
}
\end{remark}

\noindent {\it Proof of Theorem \ref{lem.F}}.
Here we give a sketch of the proof of Theorem \ref{lem.F}. 
The unique solution $u$ to \eqref{S_a} decaying at spatial infinity is expressed as $u(x) = \int_{\R^2} \Gamma_a (x,y) f (y) \dd y$, and we focus on the proof of \eqref{thelemma3} and \eqref{1.4}. By the integration by parts we have 
\begin{align*}
\int_{\R^2} \Gamma_a (x,y) f (y) \dd y & \,=\, -\int_{\R^2} \nabla_y \Gamma_a (x,y) F (y) \dd y \\
& \,=\, - \bigg ( \int_{|y|<\frac{|x|}{2}} + \int_{\frac{|x|}{2}\leq |y|} \bigg ) \nabla_y \Gamma_a (x,y) F (y) \dd y \\
& \,=\,  I (x) +  II(x)\,.
\end{align*}
The term $I$ is further decomposed as 
\begin{align*}
I (x) & \,=\, -\int_{|y| < \frac{|x|}{2}} \nabla_y L(x,y) F (y) \dd y - \int_{|y| < \frac{|x|}{2}} \nabla_y \big ( \Gamma_a (x,y) - L(x,y) \big ) F (y) \dd y\\
& \,=\, I_1 (x) + I_2 (x)\,.
\end{align*}
By the definition of $L(x,y)$ we have $-\big ( \nabla_y L(x,y)\big ) F = (F_{12}-F_{21}) \frac{x^\bot}{4\pi |x|^2}$, which implies 
\begin{align}
I_1 (x) \,=\, \int_{|y| < \frac{|x|}{2}} \big (F_{12}(y) -F_{21}(y) \big ) \dd y \frac{x^\bot}{4\pi |x|^2}\,.
\label{est.I_1}
\end{align}
As for $I_2$, when $|x| \geq 1$ we have from \eqref{est.lem.thm.linear.whole.1.1} with $m=1$,
\begin{align}\label{est.I_2}
|I_2 (x) | & \leq C \min \big\{ \frac{1}{|a| |x|^3}, \frac{1}{|x|} \big\}  \int_{|y|\leq \frac{|x|}{2}} |F(y) | \dd y + \frac{C}{|x|^2} \int_{|y|\leq \frac{|x|}{2}} |y| |F(y)| \dd y \nonumber \\
& \leq  C  \min \big\{ \frac{1}{|a| |x|^3}, \frac{1}{|x|} \big\}  \int_{|y|\leq \frac{|x|}{2}} |F(y) | \dd y + \frac{C}{(1+|x|)^{1+r}} 
\frac{1}{1-r} \| F \|_{L^\infty_{2+r}}\,,
\end{align}
where $C$ is a numerical constant independent also of $r$. 
Next we have from the direct calculation 
\begin{equation*}
|(\nabla_x K)(x,t)| \leq  C \big( t^{-\frac{3}{2}} e^{-\frac{|x|^2}{16t}} 
+  \int_{t}^{\infty} s^{-\frac{5}{2}} e^{- \frac{|x|^2}{16s}} \dd s \big)\,,
\end{equation*} 
which implies 
\begin{align*}
\int_0^\infty |(\nabla K) (O(a t) x,t )| \dd t\leq \frac{C}{|x|}\,, \qquad ~x\ne 0\,.
\end{align*}
Then by the change of the variables $y = O(a t)z$ we have 
\begin{align}\label{proof.thm.linear.whole.5}
|II(x)| & \leq  \big|  \int_{|y| \ge \frac{|x|}{2} } \nabla_{y} \Gamma_{a}(x,y) F(y) \dd y \big|  \nonumber \\
& \leq     \int_{0}^{\infty} \int_{|y| \ge \frac{|x|}{2} }  |(\nabla K)(O(a t)x - y, t)|  |F(y)| \dd y \dd t \nonumber \\
& \leq C \|  F \|_{L^\infty_{2+r}} \int_{|z| \ge \frac{|x|}{2} } \bigg( \int_{0}^{\infty}  |(\nabla K)(O(a t)(x - z), t)| \dd t  \bigg) (1+|z|)^{-2-\gamma} \dd z  \nonumber \\
& \leq  C  \| F \|_{L^\infty_{2+r}} \int_{|z| \geq  \frac{|x|}{2}} |x-z|^{-1} (1+|z|)^{-2-\gamma} \dd z \nonumber \\
& \leq  \frac{C}{(1+|x|)^{1+\gamma}} \|F\|_{L^\infty_{2+r}}\,.
\end{align}
Here $C$ is a numerical constant. 
From \eqref{est.I_1}, \eqref{est.I_2}, and \eqref{proof.thm.linear.whole.5}, we conclude \eqref{thelemma2} and \eqref{thelemma3} for $r\in [0,1)$. The proof of Theorem \ref{lem.F} is complete.

\section{Proof of linear result}\label{sec.thm.main}
In this section we prove Theorem \ref{thm.main}. Set 
\begin{align}\label{th1}
			L [f](x) \,=\, \int_{\mathbb{R}^2} \Gamma_a (x,y) f(y) \dd y\,,
\end{align}
where $\Gamma_a (x,y)$ is given by \eqref{gamma}. It is known by \cite[Lemma 3.5]{H3} that $u=L[f]$ together with $p$ defined by \eqref{bangou} is the unique weak solution to \eqref{a} decaying at spatial infinity. So we focus on the proof of the estimates \eqref{1.1} and \eqref{1.2} here. To apply Lemma \ref{lem.thm.linear.whole.1} we first divide \eqref{th1} into three parts:
\begin{align*}
			L [f](x)&\,=\, U_1(x) +U_2(x) + U_3(x) \\
			&\,=\, \bigg ( \int_{|y|<\frac{|x|}{2}} + \int_{\frac{|x|}{2}\leq |y|\leq 2 |x|} +  \int_{2|x|<|y|} \bigg )\, \Gamma_a (x,y) f(y) \dd y\,.
\end{align*} 
By Lemma \ref{lem.thm.linear.whole.1} we have 
\begin{equation}\label{th2} 
U_1(x) \,=\, \int_{|y|<\frac{|x|}{2}} y^\bot \cdot f(y) \dd y \frac{x^\bot}{4\pi |x|^2} + W_1(x) 
\end{equation}
with
\begin{align}\label{th3}
|W_1(x)| &\leq C \min \big\{ \frac{1}{|a| |x|^2}, \frac{1}{|a|^{\frac{1}{2}} |x|} \big\} \int_{|y|<\frac{|x|}{2}} |f(y)| \dd y \nonumber \\
& \qquad + C \min \big\{ \frac{1}{|a| |x|^2}, 1 \big\} \int_{|y|<\frac{|x|}{2}} |f(y)| \dd y  
+ \frac{C}{|x|^{2}} \int_{|y|<\frac{|x|}{2}} |y|^2 |f(y)| \dd y \nonumber \\
&\leq 
\begin{cases}
& \displaystyle C \big ( \frac{1}{|a|^\frac12} + 1\big ) \| f \|_{L^\infty}\,, \qquad \qquad |x|\leq 1\,, \\
& \displaystyle C \big \{ \frac{1}{|a|^\frac{1+r}{2} |x|^{1+r}} + \frac{1}{(1-r)|x|^{1+r}} \big \} \| f \|_{L^\infty_{3+r}}\,, \qquad \qquad |x|\geq 1
\end{cases} \nonumber \\
& \leq \frac{C}{(1+|x|)^{1+r}} \big ( \frac{1}{|a|^\frac{1+r}{2}} + \frac{1}{1-r} \big )\| f \|_{L^\infty_{3+r}}\,.
\end{align}
Similarly, by Remark \ref{rem.lem.thm.linear.whole.1} we have 
\begin{align}\label{th4}
			|U_3(x)| &\leq  C \int_{2|x|<|y|} \bigg (  \min \big\{ \frac{1}{|a| |y|^2}, \frac{1}{|a|^{\frac{1}{2}} |y|} \big\}  + \min \big\{ \frac{1}{|a| |y|^2}, 1 \big\}  + \frac{|x|}{|y|} \bigg ) |f(y)| \dd y \nonumber \\
& \leq 
\begin{cases}
& \displaystyle C \big ( \frac{1}{|a|^\frac12} + 1 \big ) \| f\|_{L^\infty_{3+r}}\,, \qquad \qquad |x|\leq 1\,, \\
& \displaystyle \frac{C}{|x|^{1+r}} \big ( \frac{1}{|a|^\frac{1}{2}} + 1 \big ) \| f\|_{L^\infty_{3+r}}\,, \qquad \qquad |x|\geq 1\,.
\end{cases}
\end{align}
From \eqref{th4} we have 
\begin{align}\label{th4''}
|U_3 (x)| \leq \frac{C}{(1+|x|)^{1+r}} \big ( \frac{1}{|a|^\frac{1}{2}}  + 1\big ) \| f \|_{L^\infty_{3+r}}
\end{align}
with a numerical constant $C$. 
Finally, we decompose $U_2(x)$ as  
\begin{equation}\label{th5}
			U_2(x) \,=\, U_{2,1}(x) + U_{2,2}(x)\,, \qquad U_{2,1} \,=\, U_{2,1,1} + U_{2,1,2}
\end{equation}
with
\begin{align*}
			U_{2,1,1}(x) &\,=\, \int_{\frac{|x|}{2}\leq |y|\leq 2|x|} \int^l_0 O(at)^\top \frac{1}{8\pi t} e^{-\frac{|O(at)x-y|^2}{4t}} f(y) \dd t \dd y\,, \\
			U_{2,1,2}(x) &\,=\, \int_{\frac{|x|}{2}\leq |y|\leq 2|x|} \int^\infty_l O(at)^\top  \frac{1}{8\pi t} e^{-\frac{|O(at)x-y|^2}{4t}}  f(y) \dd t \dd y\,, \\
			U_{2,2}(x) &\,=\, \int_{\frac{|x|}{2}\leq |y|\leq 2|x|} \int^\infty_0 O(at)^\top \big ( K(O(at) x-y) -  \frac{1}{8\pi t} e^{-\frac{|O(at)x-y|^2}{4t}} \mathbb{I} \big ) f(y) \dd t \dd y\,,
\end{align*}
where $l=l(a,|x|)>0$ will be chosen later.
We start from the estimate of $U_{2,1,1}(x)$.
By Fubini's theorem and changing the variable as $z=O(at)x-y$ we obtain
\begin{align}\label{th6.1}
|U_{2,1,1}(x)| & \leq \frac{C}{(1+|x|)^{3+r}} \|f\|_{L^\infty_{3+r}} \int_{\frac{|x|}{2}\leq |y|\leq 2|x|} \int^l_0 t^{-1}e^{-\frac{|O(at)x-y|^2}{4t}} \dd t \dd y \nonumber \\
			&\leq \frac{C}{(1+|x|)^{3+r}} \|f\|_{L^\infty_{3+r}} \int^l_0 \int_{\R^2} t^{-1} e^{-\frac{|O(at)x-y|^2}{4t}} \dd y  \dd t \nonumber \\
			&\leq \frac{C}{(1+|x|)^{3+r}} \|f\|_{L^\infty_{3+r}} \int_0^l  \int_{\R^2} t^{-1} e^{-\frac{|z|^2}{4t}} \dd z \dd t  \nonumber \\
			&\leq \frac{C}{(1+|x|)^{3+r}} l \, \|f\|_{L^\infty_{3+r}} \,.
\end{align}
Here $C$ is a numerical constant. Next we estimate $U_{2,1,2}$.
Since
\begin{align*}
O(at)^\top \,=\, -\frac1{a} \frac{\dd}{\dd t} \dot{O}(at)^\top\,,
\end{align*}
the integrating by parts yields
\begin{align}\label{th7.1}
U_{2,1,2}(x) &= -\frac1{2a} \int_{\frac{|x|}{2}\leq |y|\leq 2|x|} \int^\infty_l \big( \frac{\dd}{\dd t} \dot{O}(at)^\top \big ) G(O(at)x-y,t) f(y) \dd t \dd y \nonumber \\
&= \frac1{2a} \int_{\frac{|x|}{2}\leq |y|\leq 2|x|} \int^\infty_l \dot{O}(at)^\top \frac{\dd}{\dd t} \big ( G(O(at)x-y,t) \big )  f(y) \dd t \dd y + W_2(x)\,,
\end{align}
and the remainder term $W_2$ is estimated as
\begin{align}\label{th7.2}
|W_2(x)| &\leq \frac{C}{|a|} \int_{\frac{|x|}{2}\leq |y|\leq |x|} |G(O (a l) x-y, l) f(y)| \dd y \nonumber \\
& \leq \frac{C}{(1+|x|)^{1+r}} \frac{1}{l |a|} \| f\|_{L^\infty_{3+r}} \,.
\end{align}
To estimate the first term in the right-hand side of \eqref{th7.1} we use the following calculation,
\begin{align*}
& \frac{\dd}{\dd t} G(O(at)x-y,t) \\
&\quad =\frac{e^{-\frac{|O(at)x-y|^2}{4t}}}{4\pi} 
\big\{ -t^2 + t^{-3} \frac{|O(at)x-y|^2}{4} - a t^{-2} \frac{(\dot O(at)x) \cdot (O(at)x - y)}{2} \big\} \,.
\end{align*}
Hence we have 
\begin{align*}
\big | \int_{\frac{|x|}{2}\leq |y|\leq 2|x|} \int^\infty_l \dot{O}(at)^\top \frac{\dd}{\dd t} \big ( G(O(at)x-y,t) \big ) f(y) \dd t \dd y \big | \hspace{-9cm} \\
			&\leq  C \int_{\frac{|x|}{2}\leq |y|\leq 2|x|} \int^\infty_l \big ( t^{-2} + t^{-3}|O(at)x-y|^2 \big ) e^{-\frac{|O(at)x-y|^2}{4t}} |f(y)| \dd t \dd y \\
			&\quad + C \ |a| \int_{\frac{|x|}{2}\leq |y|\leq 2|x|} \int^\infty_l t^{-2} |x| |O(at)x-y| e^{-\frac{|O(at)x-y|^2}{4t}} |f(y)| \dd t \dd y \\
                   &\leq \frac{C}{(1+|x|)^{3+r}} \|f\|_{L^\infty_{3+r}} \big ( \frac{1}{l} |x|^2 +|a||x| \int_{\frac{|x|}{2}\leq |y|\leq 2|x|} \int_0^\infty t^{-\frac32} e^{-\frac{|O(at)x-y|^2}{8t}} \dd t \dd y \big )\,,
\end{align*}
and then, by the change of variables as $z = O(at)^\top y$ we see 
\begin{align}\label{th7.3}
&\leq \frac{C}{(1+|x|)^{3+r}} \|f\|_{L^\infty_{3+r}} \big (\frac1{l} |x|^2 + |a||x| \int_{\frac{|x|}{2}\leq |z|\leq 2|x|} \int_0^\infty   t^{-\frac32}e^{-\frac{|x-z|^2}{8t}} \dd t \dd z \big ) \nonumber \\
&\leq \frac{C}{(1+|x|)^{3+r}} \|f\|_{L^\infty_{3+r}} \big (\frac1{l} |x|^2 + |a||x| \int_{|x-z|\leq3|x|} \frac{\dd z}{|x-z|} \big )  \nonumber \\
&\leq \frac{C}{(1+|x|)^{1+r}} (\frac{1}{l} + |a|) \| f\|_{L^\infty_{3+r}}\,.
\end{align}
Then \eqref{th7.1}, \eqref{th7.2}, and \eqref{th7.3} implies that
\begin{equation}\label{th7}
			|U_{2,1,2}(x)| \leq \frac{C}{(1+|x|)^{1+r}} (\frac{1}{l|a|} + 1) \|f\|_{L^\infty_{3+r}}\,.
\end{equation}
Here $C$ is a numerical constant. On the other hand, the term  $U_{2,2}$ converges absolutely without using the effect of rotation. Indeed, changing the variables $y = O(at) z$, we have
\begin{align}\label{th8.1}
\begin{split}
\big | \int_{\frac{|x|}{2}\leq |y|\leq 2|x|} \int^\infty_0 O(at)^\top \big ( K(O(at)x-y,t) - \frac{1}{8\pi t} e^{-\frac{|O(at)x-y|^2}{4t}} \mathbb{I} \big ) f(y) \dd t \dd y \big | \\
			\leq \frac{C}{(1+|x|)^{3+r}} \|f\|_{L^\infty_{3+r}} \int_{|z|\leq2|x|} \int^\infty_0 |B(x-z,t)| \dd t \dd z \,, 
\end{split}
\end{align}
where $B(x,t)$ is given by 
\begin{align*}
B(x,t) \,=\, \bigg ( \frac{e^{-\frac{|x|^2}{4t}}}{8\pi t} - \frac{1-e^{-\frac{|x|^2}{4t}}}{2\pi |x|^2} \bigg ) \bigg ( \mathbb{I} - 2\frac{x\otimes x}{|x|^2} \bigg )\,.
\end{align*}
For any fixed $x-z$ we have from the change of variables as $s = \frac{|x-z|^2}{4t}$,
\begin{align}\label{th8.2}
\int_0^\infty |B(x-z,t) | \dd t \leq \big | \mathbb{I} - 2\frac{(x-z)\otimes (x-z)}{|x-z|^2} \big | \int_0^\infty \frac{-s e^{-s} + 1-e^{-s}}{s^2} \dd s \leq C\,,
\end{align}
where $C$ is independent of $x-z$. Here we have used the identity
\begin{align*}
\frac{\dd}{\dd s} \big ( \frac{e^{-s}-1}{s} \big )
\,=\, \frac{-s e^{-s} + 1-e^{-s}}{s^2} > 0\,.
\end{align*}
We combine \eqref{th8.1} with \eqref{th8.2} to conclude
\begin{align}\label{th8}
|U_{2,2}(x)| \leq \frac{C}{(1+|x|)^{3+r}} \|f\|_{L^\infty_{3+r}}\,.
\end{align}
Here $C$ is a numerical constant. Collecting \eqref{th5}, \eqref{th6.1}, \eqref{th7}, and \eqref{th8}, we see
\begin{align*}
|U_2(x)| \leq \frac{C}{(1+|x|)^{1+r}} \big\{ \frac{l}{(1+|x|)^2} + \frac{1}{l|a|} + 1 \big\} \| f\|_{L^\infty_{3+r}}\,,
\end{align*}
and thus, by taking $l = \frac{1+|x|}{|a|^{\frac12}}$,
\begin{align}\label{th9}
|U_2(x)| \leq \frac{C}{(1+|x|)^{1+r}}  \big\{ \frac{1}{(1+|x|) |a|^\frac12}  + 1 \big\} \| f\|_{L^\infty_{3+r}}\,.
\end{align}
From \eqref{th1}, \eqref{th2}, \eqref{th3}, \eqref{th4''}, and \eqref{th9}, we obtain \eqref{1.1} and \eqref{1.2}. The proof of Theorem \ref{thm.main} is complete.

\section{Proof of nonlinear result}\label{sec.nonlinear}
We are now in a position to give a proof of our main result. The unique existence and the asymptotic behavior of solutions to \eqref{NS_a} will be obtained by combining the results of Theorem \ref{lem.F} and Theorem \ref{thm.main} by applying the standard fixed point argument. 
For $r\in [0,1)$ and $\delta\in (0,1)$ we introduce the function space $X_{r,\delta}$ as follows.
\begin{align*}
X_{r,\delta} \,=\, \{ v \in L^\infty_{1+r}(\R^2)^2 ~|~ \ ||v||_{L^\infty_{1+r}}\leq\delta, \ \ {\rm div}\ v= 0 \}\,.
\end{align*}
We also set
\begin{align}\label{the1}
\begin{split}
U(x) & \,=\, \frac1{2\pi} \frac{x^\bot}{|x|^2}(1-e^{-\frac{|x|^2}4})\,,\\
\alpha & \,=\, \frac{\int_{\R^2} y^{\bot} \cdot f(y)\ \dd y}{\int_{\R^2} y^{\bot} \cdot \Delta U(y) \ \dd y} 
\,=\, \frac12 \int_{\R^2} y^{\bot} \cdot f(y)\ \dd y\,, \\
w(x) & \,=\, u(x) - \alpha U(x)\,,
\end{split} 
\end{align}
Here we have used the fact $\int_{\R^2} x^\bot \cdot \Delta U \dd x = 2$, which is derived from the identity 
$$\Delta U \,=\, (-\partial_2 G\,, \partial_1 G)^\top\,, \qquad  G(x) \,=\, \frac{1}{4\pi} e^{-\frac{|x|^2}{4}}\,. $$
The direct computation leads to the existence of a scalar function $P_U \in L^\infty (\R^2)$ such that
\begin{align*}
-a\big (x^\bot \cdot \nabla \alpha U - \alpha U^\bot \big ) + \alpha^2 U \cdot \nabla U \,=\, \nabla P_U\,.
\end{align*}
Then $w$ satisfies the following equations in $\R^2$:
\begin{align*}
\begin{cases}
		-\Delta w -a(x^\bot \cdot \nabla w-w^\bot) + \nabla\pi \\
		\qquad \,=\, -\alpha \big (U\cdot \nabla w + w\cdot \nabla U \big ) - w\cdot \nabla w -\alpha \Delta U + f\,, \\
		{\rm div}\, w \,=\, 0\,.
\end{cases}
\end{align*}
Here $\pi = p - P_U$. Let us recall that for $f \in L^\infty_{3} (\R^2)^2$, the function $L[f](x) = \int_{\R^2} \Gamma_a (x,y) f(y) \dd y$ defines the unique weak solution to \eqref{S_a} decaying at spatial infinity. Then we introduce the map $\Phi$ as
\begin{align}\label{the4}
			\Phi[w](x) \,=\, L\big[-\alpha(U\cdot \nabla w + w\cdot \nabla U) - w\cdot \nabla w -\alpha \Delta U + f\big](x)\,.
\end{align}
Here we consider the leading profile of \eqref{the4}. Since $U\cdot \nabla w + w\cdot \nabla U = {\rm div}\, (U\otimes w +w\otimes U)$ and $w\cdot \nabla w= {\rm div}\, ( w\otimes w )$, we see that
\begin{align}\label{the5}
\begin{split}
				\frac1{4\pi}&\left\{ \int_{|y|<\frac{|x|}{2}} (U\otimes w +w\otimes U)_{2,1}-(U\otimes w +w\otimes U)_{1,2}\ \dd y \right. \\
				&\qquad \left. +  \int_{|y|<\frac{|x|}{2}} (w\otimes w)_{2,1} - (w\otimes w)_{1,2}\ \dd y \right\} \\
				\,=\,0\,. 
\end{split}
\end{align}
From \eqref{the4} and \eqref{the5}, Theorems \ref{thm.main} and \ref{lem.F} yield
\begin{align}\label{phi}
\begin{split}
			\Phi[w](x) \,=\, & \mathcal{R}[-\alpha (U\cdot \nabla w + w\cdot \nabla U) - w\cdot \nabla w -\alpha \Delta U + f](x) \\ 
			&+ \bigg( \int_{|y|<\frac{|x|}{2}} y^\bot \cdot f \dd y - \alpha \int_{|y|<\frac{|x|}{2}} y^\bot \cdot \Delta U \dd y\bigg) \frac{x^\bot}{4\pi |x|^2}\,.
\end{split}
\end{align}
To estimate the last term of \eqref{phi}, we have from the definition of $\alpha$ in \eqref{the1} and $\int_{\R^2} x^\bot \cdot \Delta U \dd x =2$,
\begin{align}\label{yokei}
& \big | \int_{|y|<\frac{|x|}{2}} y^\bot \cdot f \dd y - \alpha \int_{|y|<\frac{|x|}{2}} y^\bot \cdot \Delta U \dd y \big| \nonumber \\
& \,=\, \frac12 \big| \int_{\R^2} y^\bot \cdot \Delta U \dd y \ \int_{|y|<\frac{|x|}{2}} y^\bot \cdot f \dd y 
- \int_{|y|<\frac{|x|}{2}} y^\bot \cdot \Delta U \dd y\ \int_{\R^2} y^\bot \cdot f \dd y \big| \nonumber \\
& \,=\, \frac12 \big| \int_{\R^2}y^\bot \cdot \Delta U \dd y\ 
        \bigg( 
					\int_{|y|<\frac{|x|}{2}} y^\bot \cdot f \dd y - \int_{\R^2} y^\bot \cdot f \dd y 
				\bigg) \nonumber \\
			&\qquad-  \bigg(
			 		\int_{|y|<\frac{|x|}{2}} y^\bot \cdot \Delta U \dd y - \int_{\R^2} y^\bot \cdot \Delta U \dd y 
			 	\bigg)
			 \ \int_{\R^2} y^\bot \cdot f \dd y \big| \nonumber \\
			& \,=\, \frac12 \big|-\int_{\R^2} y^\bot \cdot \Delta U \dd y \int_{\frac{|x|}2\leq|y|} y^\bot \cdot f \dd y + \int_{\frac{|x|}2\leq|y|} y^\bot \cdot \Delta U \dd y \int_{\R^2} y^\bot \cdot f \dd y \big| \nonumber \\ 
			&\leq \frac{C_r}{(1+|x|)^{r}} (\|y^\bot\cdot f\|_{L^1} +\|f\|_{L^\infty_{3+r}} )\,, \qquad {\rm for}~~ |x|>1\,.
\end{align}
Here the constant $C_r$ depends only on $r$. Note that the boundedness of $\|y^\bot\cdot f\|_{L^1}$ is always valid when $r>0$. When $|x|\leq 1$ it is easy to see 
\begin{align}\label{yokei'}
\big| \int_{|y|<\frac{|x|}{2}} y^\bot \cdot f \dd y - \alpha \int_{|y|<\frac{|x|}{2}} y^\bot \cdot \Delta U \dd y \big|   \leq C (\| y^\bot\cdot f\|_{L^1} + \|f\|_{L^\infty} )\,.
\end{align}
Estimates \eqref{yokei} and \eqref{yokei'} combined with \eqref{phi} imply that
\begin{align}\label{nokori}
\begin{split}
		\Phi[w](x) \,=\, & \mathcal{R}[-\alpha(U\cdot \nabla w + w\cdot \nabla U) - w\cdot \nabla w -\alpha \Delta U + f](x) \\ 
		& + O\bigg ((  \|y^\bot\cdot f\|_{L^1} +\|f\|_{L^\infty_{3+r}})(1+|x|)^{-1-r} \bigg )\,.
\end{split}
\end{align}
By \eqref{nokori} and Theorems \ref{thm.main} and \ref{lem.F}, we can verify that
\begin{align}\label{saigo}
\|\Phi[w] \|_{L^\infty_{1+r}} &\leq C \big (\frac{1}{1-r} + \frac{1+\log |a|}{|a|^\frac{r}{2}} \big ) \| \alpha (U \otimes w + w \otimes U) + \frac12 w\otimes w \|_{L^\infty_{2+r}} \nonumber \\
& \quad + C \big (\frac{1}{1-r} + \frac{1}{|a|^\frac{1+r}{2}} \big ) \|-\alpha \Delta U + f \|_{L^\infty_{3+r}} \nonumber \\
& \qquad + C_r  \big ( \|y^\bot\cdot f\|_{L^1} +\|f\|_{L^\infty_{3+r}} \big) \nonumber \\
\begin{split}
&\leq C \big (\frac{1}{1-r} + \frac{1+\log |a|}{|a|^\frac{r}{2}} \big ) \big ( |\alpha| \|w\|_{L^\infty_{1 + r}}+ \|w\|_{L^\infty_{1 + r}}^2 \big ) \\
& \quad + C_r \big ( \frac{1}{1-r} + \frac{1}{|a|^\frac{1+r}{2}} \big ) \big ( \| f \|_{L^\infty_{3+r}} + \|y^\bot\cdot f\|_{L^1} \big )\,,
\end{split} 
\end{align}
where $C_r$ depends only on $r$, while $C$ is a numerical constant. We may take $C$ and $C_r$ larger than $1$, and note that $|\alpha|\leq 2^{-1}\| y^\bot \cdot f \|_{L^1}$. If 
$$0<\delta \leq \frac{1}{3C (\frac{1}{1-r}+\frac{1+\log |a|}{|a|^\frac{r}{2}})}$$
and 
\begin{equation*}
\lambda(f) \,=\,  \| f\|_{L^\infty_{3+r}} + \|y^\bot\cdot f\|_{L^1} \leq \frac{\delta}{3C_r (\frac{1}{1-r}+\frac{1}{|a|^\frac{1+r}{2}})}\,,
\end{equation*}
then we see that $\Phi[w]$ becomes a mapping from $X_{r,\delta}$ into $X_{r,\delta}$. 
Moreover, from \eqref{1.4} there is a numerical constant $C'>0$ such that 
\begin{align*}
& \|\Phi[w_1] - \Phi[w_2]\|_{L^\infty_{1 + r}} \\
& \quad \,=\,  \| \mathcal{R} \big[-\alpha \big\{ U\cdot \nabla (w_1 - w_2) + (w_1 - w_2)\cdot \nabla U \big\} 
 - w_1\cdot \nabla w_1 +w_2\cdot \nabla w_2 \big] \|_{L^\infty_{1+ r}} \\
& \quad \leq C' \big (\frac{1}{1-r}+\frac{1+\log |a|}{|a|^\frac{r}{2}} \big ) \|
				-\alpha\big\{ (U\otimes (w_1 - w_2) + (w_1 - w_2)\otimes U) \big\} \\
& \qquad \qquad \qquad \qquad \qquad \qquad  \quad
- \frac12 (w_1 \otimes w_1 - w_2 \otimes w_2) \|_{L^\infty_{2+r}} \\
& \quad \leq C' \big (\frac{1}{1-r}+\frac{1+\log |a|}{|a|^\frac{r}{2}} \big ) 
(|\alpha| + \|w_1\|_{L^\infty_{1+ r}} + \|w_2\|_{L^\infty_{1+ r}}) 
\|w_1 - w_2\|_{L^\infty_{1 + r}}  \\
& \quad \leq C' \big (\frac{1}{1-r} +\frac{1+\log |a|}{|a|^\frac{r}{2}} \big )\big(\frac{\lambda(f)}{2} +2\delta \big)  \|w_1 - w_2\|_{L^\infty_{1 + r}}  \\
& \quad \leq 3\delta C' \big (\frac{1}{1-r} +\frac{1+\log |a|}{|a|^\frac{r}{2}}\big )  \|w_1 - w_2\|_{L^\infty_{1 + r}}  \\
& \quad \,=\, \tau  \|w_1 - w_2\|_{L^\infty_{1 + r}}\,,
\end{align*}
for all $w_1,w_2 \in X_{r,\delta}$, where we have set $\tau = 3 \delta C'(\frac{1}{1-r}+\frac{1+\log |a|}{|a|^\frac{r}{2}})$. 
Hence, if  $\delta$ (and thus, also $\lambda(f)$) is sufficiently small so that $\tau\in (0,1)$ is justified, then we can conclude that $\Phi$ is a contraction on $X_{r,\delta}$.
By the fixed point theorem, there exists a fixed point $v$, which is unique in $X_{r,\delta}$, such that
\begin{align*}
u(x) \,=\, \alpha U(x) + v(x)\,, \qquad v \in X_{r,\delta}
\end{align*}
is a unique solution to \eqref{NS_a} with the pressure $p$ defined by \eqref{bangou2}.
Finally, the estimate \eqref{1.5} follows from \eqref{saigo} for the fixed point $v$ of $\Phi$ by virtue of the smallness of $|\alpha|$ and $\delta$. The proof of Theorem \ref{thm.main2} is complete.

\end{document}